\documentclass[mathcal]{article}
\usepackage{amsmath}
\usepackage{amsfonts}
\usepackage{amssymb}
\usepackage{latexsym}
\usepackage{euscript}
\usepackage{diagrams}

\def\lra{\longrightarrow}
\def\frak{\mathfrak}
\def\wt{\widetilde}
\def\wtp{\wt{\rho}}
\def\text{\textrm}
\def\K{K}
\def\M{\frak M}

\def\ord{\mathrm{ord}}
\def\GL{\mathrm{GL}}
\def\IQ{I}
\def\Spec{\text{Spec}}

\def\Bbb{\mathbb}
\def\Q{\Bbb Q}
\def\Z{\Bbb Z}
\def\l{\ell}
\def\Qbar{\overline{\Q}}
\def\Kbar{\overline{K}}
\def\Gal{\mathrm{Gal}}

\def\Rem{\noindent \bf Remark\rm. }
\def\Cau{\noindent \bf Caution\rm. }
\def\Exa{\noindent \bf Example\rm. }
\def\Proof{\noindent \bf Proof\rm. }
\def\mod{\text{ mod }}

\def\qed{\square}

\begin{document}
\title{Almost Rational Torsion Points on  Semistable Elliptic Curves}
\author{Frank Calegari}

\maketitle
\abstract{If $P$ is an algebraic point on a commutative group scheme $A/K$,
then $P$ is \emph{almost rational} if no two
non-trivial
Galois conjugates $\sigma P$, $\tau P$ of $P$ have sum  equal to
$2P$.
In this paper, we classify almost rational torsion points on
semistable elliptic curves over $\Q$.}
\newtheorem{theorem}{Theorem}[section]
\newtheorem{lemma}{Lemma}[section]
\newtheorem{cor}{Corollary}[section]
\newtheorem{df}{Definition}
\section{Definitions and Results}

Let $X$ be an algebraic curve of genus greater than one. Let
$J(X)$ be the Jacobian variety of $X$, and embed $X$ in $J(X)$.
The Manin-Mumford conjecture states that
the set
of torsion points
$X_{\mathrm{tors}}:=X \cap J_{\mathrm{tors}}$ is finite. This conjecture
was first proved in 1983 by Raynaud \cite{Raynaud}. It has long been 
known that the geometry of $X$ imposes  strong 
conditions on the action of Galois on $X_{\mathrm{tors}}$.
An approach to the Manin-Mumford conjecture using 
Galois representations attached to Jacobians
was first suggested
by Lang \cite{Lang}.
Recently, by exploiting the relationship 
between 
the action of
Galois on modular Jacobians and
the Eisenstein ideal 			
(developed by Mazur \cite{Mazur3})
Baker \cite{Baker} and (independently) 			
Tamagawa \cite{Tamagawa}
explicitly determined the set of torsion points of $X_0(N)$ for
$N$ prime.
Developing  these ideas further,
Ribet defined
the notion of an \emph{almost rational torsion point},
and used this concept to derive the
Manin-Mumford conjecture \cite{Ribet4}
using some unpublished results of Serre \cite{SerreU}.  
One idea suggested by these papers is that a possible approach
to finding all torsion points on a curve $X$ is to determine the
set of
almost rational torsion points on $J(X)$. 
Moreover, the concept of an almost rational torsion point makes sense
for any Abelian variety, or more generally any commutative group scheme, and
the set of such points may be interesting to study in their own right.
In this paper we consider almost rational points on semistable elliptic
curves over $\Q$, and prove (in particular) that they are all defined
over $\Q(\sqrt{-3})$.

\newpage

Let $K$ be a perfect field, and let
$A/K$ be a commutative algebraic group scheme.
\begin{df} Let $P \in A(\Kbar)$. Then
$P$ is \emph{almost rational over $K$}  if the
following  holds:
For all   $\sigma$, $\tau \in \Gal(\Kbar/K)$,  the equality
$$\sigma P + \tau P = 2 P$$
is satisfied only when  $P = \sigma P = \tau P$.
Equivalently, $P$ is almost rational if any two non-trivial
Galois conjugates of $P$ have sum different from $2P$.
If $P$ is also a torsion point of $A$, then we call
$P$ an \emph{almost rational torsion point}. 
\end{df}
\begin{lemma} Let $P$ be an almost rational point over $K$.
\begin{enumerate}
\item If $L \supset K$ is a subfield of $\Kbar$, then $P$ is almost
rational over $L$.
\item If $Q$ is a $K$-rational point of $A$,
 then $Q$ is almost rational over $K$.
\item If $Q$  is a $K$-rational point of $A$,
 then $P + Q$ is almost rational over $K$.
\item Let $\gamma \in \Gal(\Kbar/K)$. Then $\gamma P$ is almost
rational over $K$.
\item Let $\sigma \in \Gal(\Kbar/K)$. 
If $(\sigma - 1)^2 P = 0$ then $\sigma P = P$.
\item If $2P$ is $K$-rational, then so is $P$.
\end{enumerate}
\end{lemma}
\Proof The statements $1$, $2$ and $3$ are obvious.  For $4$, suppose that
$\sigma \gamma P + \tau \gamma P = 2 \gamma P$.
Then
$$\gamma^{-1} \sigma \gamma P + \gamma^{-1} \tau \gamma P = 2P.$$
Since $P$ is almost rational,
 $P = \gamma^{-1} \sigma \gamma P = \gamma^{-1} \tau \gamma P$
and thus $\gamma P = \sigma \gamma P = \tau \gamma P$. 
For 5, one has $\sigma^2 P - 2 \sigma P + P = 0$. Applying 
$\sigma^{-1}$ to the left hand side of this equation, one finds that
$\sigma P + \sigma^{-1} P = 2P$.
Setting $\tau = \sigma^{-1}$ and using
the fact that $P$ is almost rational we conclude that $P = \sigma P$.
For 6, let $\sigma P$ be any Galois conjugate of $P$. Then
$$\sigma P + \sigma P = 2 \sigma P = \sigma(2P) = 2P$$
and so by almost rationality, $P = \sigma P$ for all
$\sigma$, and so $P$ is $K$-rational.
$\qed$

\Cau We shall see below (Theorem $1.2$) that the set of
almost rational points does not necessarily form a group. In fact, the
almost rationality of $P$ does not imply the almost rationality of
multiples of $P$.

\

\begin{lemma} Let $A = \Bbb G_m/\Q$. Then the almost rational torsion
points on $A$ are exactly the points of order dividing $6$. Let
$H = \mu_n/K$, where $K$ is any field such that
$\Gal(K(\zeta_n)/K) \simeq \Gal(\Q(\zeta_n)/\Q)$. 
 Then the almost rational torsion points on
$H$ are the points of order dividing $6$. 
\end{lemma}

\Proof Let $P$ be a torsion point of order $n$ on $A$.
The Galois module generated by $P$ is the cyclotomic
module $\mu_n$.
The action of Galois on $\mu_n$ is via the mod $n$
cyclotomic character: $\chi_n: \Gal(\Qbar/\Q) \rightarrow (\Z/n \Z)^*$; 
i.e.,
we have
$\sigma P = \chi_n(\sigma) P$ for $\sigma \in \Gal(\Qbar/\Q)$. 
 To show that $P$
is \emph{not}
almost rational,
it suffices to find $\sigma$, $\tau$ such
that $\chi_n(\sigma) P + \chi_n(\tau) P = 2P$ and $P \ne \sigma P$.
In particular, since $\chi_n$ acts faithfully on $\mu_n$, it
suffices to find $a$, $b \in (\Z/n \Z)^*$ such that $a \ne 1$ and
$a + b = 2$. We attempt to do this using the Chinese Remainder Theorem.
If $p^k \| n$ and $p \ne 3$ let $a \equiv 3
\mod p^k$ and $b \equiv -1 \mod p^k$. If $3^k \| n$ let
$a \equiv 4 \mod 3^k$ and $b \equiv -2 \mod 3^k$. As long as
$n \nmid 6$, we find that $a \ne 1$, and so $P$ fails to be
almost rational. If $P$ is of order $1$ or $2$, then $P$ is
rational and so almost rational by Lemma $1.1$(2). If $P$ is of
order $3$ (respectively $6$), then the only non-trivial Galois conjugate
of $P$ is $\sigma P = 2P$ (resp. $\sigma P = 5P$). In either
of these cases $\sigma P + \sigma P \ne 2P$, and so $P$
is almost rational.  An similar argument proves the result for $H$.
$\qed$  

\

The following lemma provides a connection between
almost rational points 
and torsion points on curves, as well as
providing
a natural source of almost rational points.

\begin{lemma}
Let $X/K$ be a curve of genus $g \ge 2$, and
$Q \in X(K)$  a $K$-rational
point on $X$. Let $J/K$ be the Jacobian of $X$, and
let $i_Q:X \hookrightarrow J$
be the Albanese map $P \mapsto [P - Q]$ defined over $K$. Then
for any point $P \in X(\Kbar)$, either $P$ is a hyperelliptic
branch point of $X$, or
$i_Q(P)$ is almost rational over $K$.
\end{lemma}
\Proof Since $i_Q$ is a closed immersion, we shall
identity points of $X$ with
their images. Assume that $\sigma P  + \tau P = 2 P$. Then
since $J(\Kbar) \simeq \mathrm{Pic}^{0}(X(\Kbar))$,
the divisor $D = (\sigma P) + (\tau P) - 2(P)$ is principal, and
so equals $(f)$ for some function $f$. Either $f$ is constant, in
which case $P = \sigma P = \tau P$, or $f$ is of
degree $2$, in which case $P$ is a hyperelliptic branch point of $f$.
The result follows. $\qed$

\

One of the main motivations for studying almost rational torsion
points is  the following result of Ribet \cite{Ribet4}:

\begin{theorem} Let $A/K$ be an Abelian variety. Then the set
of almost rational torsion points on $A$ is finite.
\end{theorem}
\Rem 
Ribet's proof depends on some of Serre's  unpublished
``big-image'' Theorems \cite{SerreU}
and is not effective. Note that the Manin-Mumford conjecture
follows from Theorem $1.1$ when $X$ is embedded in its
Jacobian via an Albanese map.
 This paper provides an effective
version of Theorem $1.1$ when $A/\Q$ is a
semistable elliptic curve, and in this case we give a
complete classification of  the
possible almost rational torsion points which can arise.

\

In contrast, the following result (pointed out by Bjorn
Poonen) shows that there is an abundance of almost rational
points that are \emph{not} torsion points.

\begin{lemma} For any $P \in A(\Qbar)$, there exists an $n > 0$ such
that $n P$ is almost rational.
\end{lemma}
\Proof Since $P$ is algebraic, it is defined over some Galois 
field $K/\Q$, and
has only finitely many Galois conjugates.
Choose $n$ such that every torsion point of the form $\sigma P - \tau P$
is killed by $n$.  Then $Q = n P$ is almost rational.
Indeed, by construction, the Galois orbit $S$ of $Q$ injects into
$A(K) \otimes \Bbb R$. If $Q'$ is some extremal point of the
convex hull of $S$ (with respect to the canonical height) then
the identity $\sigma Q' + \tau Q' = 2 Q'$ can hold only if
$Q' = \sigma Q' = \tau Q'$. Thus $Q'$ is almost rational. Since
$Q'$ is a Galois conjugate of $Q$, it follows that $Q$ is also
almost rational, by Lemma $1.1$(4). $\qed$

\

The main result of this paper is the following:

\begin{theorem} Let $P$ be an almost rational torsion point
over $\Q$ on a semistable elliptic curve $E/\Q$.
Then either $P$ is rational, or it can be written as a sum
$Q + R + S$, where
\begin{enumerate}
\item \ $Q$ generates a (cyclotomic) $\mu_3$ subgroup of $E[3]$.
\item \ $R$ is a rational point of $3$-power order: $R \in E(\Q)[9]$, 
\item \ $S$ is a $\Q(\zeta_3)$ rational point of $2$-power order:
$S \in E(\Q(\sqrt{-3}))[16]$.
\end{enumerate}
Moreover, all sums of this form are almost rational.
\end{theorem}
\begin{Exa} Let $E$ be the elliptic curve:
$$y^2 + yx + y = x^3 + 354 x + 4684$$
of conductor $N=1302=2\cdot3\cdot7\cdot31$.  Then as Galois modules,
$E[3] \simeq
\mu_3\oplus\Z/3\Z$, and the point $S = (-3 -6\sqrt{-3}, 5 + 12
\sqrt{-3})$ is a torsion point of order $2$.  If $Q$ is a
$3$-torsion point that generates $\mu_3$, then $P = Q+S$ is an
almost rational torsion point on $E$. 
 Note that $3P = S$ is \emph{not} an
almost rational torsion point.  
\end{Exa} 

\

The key idea of the proof of Theorem $1.2$ is to 
limit the ramification
of $E[\l]$ for the largest prime divisor $\l$ of $|P|$. 
Ribet's level lowering result \cite{Ribet1} can then 
be used to show that $E[\l]$ is reducible. There are several
reasons why we restrict our attention to elliptic curves instead of
higher dimensional Abelian varieties.
In particular, we use many strong results about representations
arising from elliptic curves (such as modularity) which have
no convenient equivalent in higher dimensions. Secondly, we rely
on the explicit determination, by Mazur \cite{Mazur3}, of the
possible rational torsion subgroups in $E(\Q)$. Some results, however,
do apply in more generality, such as Lemma $1.5$ below.

\begin{df} Let $P$ be an almost rational torsion point. Then
the Galois module $M$ associated to $P$ is the module generated
by $P$ and all its Galois conjugates:
$$M = \sum \Z \cdot \sigma P$$
Since $P$ is a torsion point, $M$ is finite as a Galois module.
\end{df}
\Exa If $P$ is rational of order $n$ then (as Galois modules)
$M \simeq \Z/n \Z$.

\Rem Throughout this paper, the word finite is used in two
different senses. `$M$ is finite as a Galois module' means
that as an abelian group, $M$ is finite. `$M$ is finite as a group
scheme' means that there exists some finite flat group scheme
$\M/\Spec \ \Z$ such that $\M(\Q) \simeq M$. Unless explicitly stated,
we shall reserve finite to mean finite as a group scheme. 
`$M$ is finite at a prime $p$', means that there exists some finite
flat group scheme $\M_p/\Spec \ \Z_p$ such that $\M_p(\Q_p) \simeq M$
as $\Gal(\Qbar_p/\Q_p)$ modules.
$M$ is finite if and only if it is finite at all primes $p$.
If the cardinality of $M$ is coprime to $p$, then $M$ is finite at
$p$ if and only if $M$ is unramified at $p$.
\begin{lemma} Let $A/K$ be a semistable Abelian variety, and 
let $P \in A(\Kbar)$
be an almost rational point on $A$ of order $n$. Then the Galois module $M$
generated
by $P$ is finite at all primes not dividing $n$. In other words, $M$ is
\emph{unramified} outside primes dividing $n$.
\end{lemma}
\Proof 
Let  $\frak{p}$ be a prime coprime to $n$.
Let $I_\frak{p} \subseteq \Gal(\Kbar/K)$ be a choice of inertia group
above $\frak{p}$. Then since $A$ is \emph{semistable}, the
action of $I_\frak{p}$ on the Tate module $T_{\l}(A)$ for $\l|n$ 
(or for any $(\l,\frak{p}) = 1$)
satisfies 
$(\sigma - 1)^2 = 0$, for all $\sigma \in I_{\frak{p}}$
 (see \cite{Groth}, expos\'{e} IX, Proposition 3.5).
In particular,
writing $P$ in terms of its $\l$-primary components we find that
$(\sigma - 1)^2 P = 0$. Since $P$ is almost rational, by Lemma $1.1$(5),
$\sigma P = P$. Applying the same argument to all the conjugates of
$P$ (which are still almost rational, by Lemma $1.1$(4)), we find that
$M$ is unramified  at $\frak{p}$, and we are done. $\qed$ 

\

\Rem For semistable elliptic curves, 
the fact that inertia at $p$ satisfies
$(\sigma - 1)^2 | T_{\l}(E) = 0$ for $(\l,p)=1$ can
be proved directly
without appeal to results of Grothendieck. In particular,
either $E$ has good reduction
at $p$, in which case $T_{\l}(E)$ is unramified at $p$ (by the
criterion of N\'{e}ron-Ogg-Shafarevich), or $E$ has  multiplicative
reduction at $p$, in which case the result follows from the
explicit description of Tate curves given (for example) in
(\cite{AOEC2}, Chapter V).

\section{Finiteness of $\wtp$}

Let $p$ be prime. Let $\rho_p$ denote the Galois representation associated to
the Tate module $T_p(E)$. Let $\wtp_p$ denote the mod $p$ 
representation arising from the action of Galois on $E[p]$.
One sees that $\wtp_p$ is the reduction of $\rho_p$ mod $p$. 

\

Fix a semistable elliptic curve $E/\Q$,
an almost rational torsion point $P$ of order $n$, and its associated
Galois module $M$. 
Lemma $1.5$  shows how one can control the ramification of
$M$ at primes 
away from $n$. In this section, we show how it is also somtimes
possible to control the ramification at primes $p$ dividing $n$.

\begin{theorem} Let
$p | n$ be a prime such that $E[p]$ is irreducible.
Then $\wtp_{p}$ is finite at $p$.
\end{theorem}
\Proof If $E$ has good reduction at $p$ then we are done, so we may 
assume that $E$ has multiplicative reduction at $p$.  We shall use the
criterion, due to Serre (\cite{Serre2}, Proposition 5, p. 191),
 that $\wtp_p$ is finite at $p$
(\emph{peu ramifi\'{e}})
if and only if $v_p(\Delta) \equiv 0 \mod p$.

\

Let $n = p^k m$ with $(m,p) = 1$.  Write $P = P_p + P'$, where
$P_p$ is of order $p^k$, and $P'$ of order $m$.
Since $M \cap E[p] \ne 0$
and $E[p]$ is irreducible, $E[p] \subseteq M$. 
Let $E_q$ denote the Tate curve isomorphic to
$E$ over $\K$, where $[\K:\Q_p]$ is an unramified extension of
degree $\le 2$. 
Since 
\mbox{$\K(\zeta_{n})(E_q[n]) = \K(\zeta_{n}
,q^{1/n})$} we have the following  diagram of fields:
$$
\begin{diagram}
\K(E_q[n]) = &  \K(\zeta_{n},q^{1/n}) & & &\\
&  \dLine & \rdLine & & \\
 &  & & \K(\zeta_{n},q^{1/m}) & \supseteq \K(E_q[m]) \\
 &  & \ldLine & & \\
 & \K(\zeta_n) & & & \\
\end{diagram}
$$
We shall show that the extension $\K(\zeta_{n},q^{1/n})/
\K(\zeta_{n},q^{1/m})$ is trivial. Assume otherwise.
Then since $\K(E_q[n])/\K$ is Galois, 
there exists a  non-trivial
Galois automorphism $\sigma \in \Gal(\K(E_q[n])/\K)$
fixing $\K(\zeta_n,q^{1/m})$. If $\sigma$ exists, then
we may take it to be of order $p$ and of the form:
$\sigma: q^{m/n} \lra  q^{m/n} \zeta_p$.
We choose a basis  $\{q^{m/n},\zeta_{p^k}\}$ 
for $E_q[p^k]$ such that $\sigma$ acts via the matrix
$$\left( \begin{array}{cc} 1 & p^{k-1} \\ 0 & 1 \end{array} \right)
\mod p^k$$ 
Since $\sigma$ fixes $E_q[m]$, $(\sigma - 1)^2 P = 0$. Considering
$\sigma$ as an element of $\Gal(\Qbar/\Q)$ via the inclusion 
$$\Gal(\Qbar_p/\K) \hookrightarrow \Gal(\Qbar/\Q)$$
we conclude from the almost rationality of $P$ that
$\sigma P = P$. Thus with respect to our chosen basis 
for $E_q[p^k]$, $P_p$
is of the form
$$\left( \begin{array}{c} a \\ bp \\ \end{array} \right)$$
From Lemma $1.1$(4) we may apply the same argument above to 
$\gamma P$ 
for every \mbox{$\gamma \in \Gal(\Qbar/\Q)$}, and so
all Galois conjugates of $P_p$ are also of this
form.
  It follows
that the mod $p$
 Galois representation $\wtp_{p}$ is upper triangular, contradicting
the irreducibility of $E[p]$. Thus
$\K(\zeta_n,q^{1/n})/\K(\zeta_n,q^{1/m})$ is the trivial extension.

It follows that $q^{1/n} \in L:=K(\zeta_{n},q^{1/m})$, and
so in particular
$q^{1/m}$ is a $p^{k}$-th power in $L$,
and thus its normalized valuation
with respect to this field satisfies $v_L(q) \equiv 0 \mod p^k$.
Since $(m,p) = 1$, the wild ramification of 
$\K(\zeta_{n},q^{1/m})$ is limited to $\K(\zeta_{p^k})$,
for which $p^{k-1} \| e(\K(\zeta_{p^k})/\Q_p)$.  Since 
$v_{L}(q) = e(L/\Q_p) \cdot v_{\Q_p}(q)$,
the valuation of $q$ with respect to $\Q_p$ must be divisible
by $p$. Thus 
$$\ord_{p}(\Delta_q) = \ord_{p}\left(q \prod_{i=1}^{\infty}
(1 - q^n)^{24} \right) = \ord_{p}(q) \equiv 0 \mod p$$
and therefore $E[p]$ is finite at $p$. $\qed$

\begin{theorem} Let $\ell$ be the largest prime factor of $n = |P|$.
Then $E[\ell]$ is reducible.
\end{theorem}

\Proof Assume otherwise.  
From Theorem $2.1$, since $E[\l]$ is irreducible,
$\wtp_{\ell}$ is finite at $\l$. 
Let $p \ne \l$ be prime. We will show that $E[\l]$ is unramified at $p$,
or that  $\l = 3$ and $p = 2$. In fact, with no assumptions on
$E[\l]$ we shall show more generally
that $M[\l^{\infty}]$ is unramified at $p$. The irreducibility
assumption on $E[\l]$ then
ensures that $E[\l] \subseteq M[\l^{\infty}]$ and thus that $E[\l]$
is  also unramified at
$p$.

For primes $p \nmid n$ the result follows from Lemma $1.5$.
Hence it suffices to consider $p \ne \l$ such that $p | n$.
Write $P$ in terms of its $q$-primary
components
$$P = P_2 + P_3 + \ldots + P_p + \ldots +   P_{\l}$$
where $P_{q}$ is a torsion point of order $q^k \| n$.
Let $\IQ_p$ denote a choice of inertia group at $p$. 
Let $\sigma \in \IQ_p$ be an element of the inertia group.
Then by the proof of  Lemma $1.5$ (or the subsequent remark),
$(\sigma - 1)^2 P_q = 0$ for $q \ne p$.
Consider the following
diagram of fields:
$$\begin{diagram} 
\Q(M[\l^{\infty}]) 	&	\rLine	&	\Q(M)	\\
\dLine		&		&	\dLine			\\
\Q		&	\rLine	&	\Q(M[p^{\infty}])	\\
\end{diagram}
$$
For $\sigma \in \IQ_p$ fixing $M[p^{\infty}]$, $\sigma P_p = P_p$ and so
$(\sigma - 1)^2 P_p = 0$.
For such $\sigma$, $(\sigma -1)^2 P = 0$, and thus
by Lemma $1.1$(5), $\sigma P = P$. Applying this to the Galois
conjugates of $P$, one concludes that
the extension 
$\Q(M)/\Q(M[p^{\infty}])$ 
is unramified at all primes above $p$. Comparing ramification indices
at $p$ in the diagram above,
$$ e_p(\Q(M)/\Q(M[\l^{\infty}]) \cdot e_p(\Q(M[\l^{\infty}])/\Q) =
e_p(\Q(M[p^{\infty}])/\Q).$$
If $P_p$ is of order $p^k$, then $M[p^{\infty}] \subseteq E[p^k]$ and
so
$e_p(\Q(M[p^{\infty}])/\Q)$ divides the order of
$\GL_2(\Z/p^k \Z) = 
(p^2-1)(p-1) p^{4k-3}$, and so 
in particular, $e:=e_p(\Q(M[\l^{\infty}]/\Q)$ also divides this number.
Yet the action of $\IQ_p$ on $T_{\l}(E)$ is unipotent, and so
$e$ is either $1$ or some power of $\l$. Since $\l$ is the
\emph{largest} prime factor of $n$, $p < \l$, and thus 
$$(\l,(p^2-1)(p-1) p^{4k-3}) = (\l,p+1)$$
which equals $1$ unless $(\l,p) = (3,2)$.
Thus if $\l \ne 3$, $M[\l^{\infty}]$ is unramified outside $\l$, 
and for all $\l$, $M[\l^{\infty}]$ is unramified outside $2$ and $\l$.
If $E[\l]$ is irreducible, then 
the Serre conductor
\cite{Serre2}  $N(\wtp_{\ell})$ is equal to $1$ or $2$ (Since $E$ is
semistable, the exponent of $2$ in the 
Serre conductor must be at most $1$).
By Wiles \cite{Wiles} and Taylor-Wiles \cite{Taylor}, $\wtp_{\ell}$ is
modular. 
Since $\wtp_{\ell}$ is finite at $\l$, For $\l > 2$
Ribet's level lowering result
\cite{Ribet1} implies that $\wtp_{\ell}$ arises from a weight  two
modular form of level one or two. 
No such form exists.  This is a contradiction, and the theorem follows.
For $\l = 2$ and $E[2]$ irreducible one finds that $N(\wtp_{2}) = 1$, 
contradicting a theorem of Tate \cite{Tate}.
$\qed$

\

\begin{cor} Let $\l$ be any prime divisor of $n$. Then $M[\l^{\infty}]$
is at most ramified at $\l$ and primes $p|n$ such that
$\l | (p^2-1)$. Moreover, if $p \ne \l$ and
 $E[p]$ is reducible, then $M[\l^{\infty}]$ can only be ramified at
$p$ if $\l |(p - 1)$.
\end{cor}

\Proof  The first part of the Corollary is proved during the
proof of Theorem $2.2$. If $E[p]$ is reducible, then
$[\Q(E[p]):\Q]$ divides $p(p-1)$, and so arguing as in Theorem
$2.2$ we conclude that $\l$ divides $(p-1)$.

\begin{cor}  Let $P$ be an almost rational torsion point of order $n$. 
Then if $\l$ is the largest prime dividing $n$, then $\l \le 7$.
\end{cor}
\Proof If $E$ is reducible at $p$, then since $E$ is semistable, it
follows from Serre 
(\cite{Serre1} Proposition 21, p. 306, and subsequent remarks) that
 either $E$ or some isogenous curve $E'$ has a rational
point of order $p$.  The  result then follows from the following
theorem of  Mazur (\cite{Mazur3}, Theorem 8, p. 35):

\begin{theorem}[Mazur] Let  $\Phi$ be the torsion subgroup of
the Mordell-Weil group of an elliptic curve over $\Q$. Then
$\Phi$ is isomorphic to one of the following $15$ groups:
$\Z/m \Z$ for $m \le 10$ or $12$; 
$\Z/2 \Z \oplus \Z/2k \Z$ for $k \le 4$.
\end{theorem}

\section{The possible cases}

Corollary $2.2$ severely limits the possible prime divisors of $n$.
In this section, we eliminate all possibilities not allowed by 
Theorem $1.2$. Let $S$ denote the set of primes dividing $n$.
Corollary $2.2$ restricts $S$ to $15$ non-trivial possibilities,
which divide into three cases.

\subsection{Case I: $S = \{p\}$, $n = p^k$} 

When $n = p^k$ our strategy is as follows. First we show 
that $M$ is
an extension of a trivial module by a cyclotomic module.
If $p \ge 3$, this sequence splits, and it suffices to
find almost rational torsion points on the cyclotomic
module $\mu_{p^r}$, which we calculated in Lemma $1.2$.
 For $p = 2$, some complications arise,
but the  essential ideas remain the same.
Some of our arguments can be shortened
using results of Tamagawa (\cite{Tamagawa}, for example Theorem $3.2$),
however,  we  proceed directly in order to be self
contained.

\
 
We begin by recalling  an elementary result from the
theory of cyclotomic fields. 
\begin{lemma} Let $h(K)$ denote the class number of the field $K$.
Suppose that $p$ is inert in $K$.
Then for $k \ge 1$,
$$p | h(K(\zeta_p)) \Leftrightarrow p | h(K(\zeta_{p^k})).$$
In words: $p$ divides $h(K(\zeta_{p^k}))$ for all $k$ if and only
if $p$ is an irregular prime with respect to $K$. 
\end{lemma}
For a proof, see (for example) Theorem $10.4$ of Washington \cite{Wash}.

\begin{lemma}
If $p \le 7$, $p \nmid h(\Q(\zeta_{p^k}))$. Moreover, $2 \nmid h(\Q(\zeta_{3 \cdot 2^k}))$ and  $2 \nmid h(\Q(\zeta_{5 \cdot 2^k}))$. For each of these
fields $K$, if $L/K$ is a Galois extension such that $\Gal(L/K)$
is a $p$-group and
$L/K$ is unramified, then $L = K$.
\end{lemma}
The first statement is a consequence of Lemma $3.1$, the fact that
the fields $\Q(\zeta_n)$ have class number $1$ for $n = 1,3,4,5,7,12,20$, and the fact that 
that $2$ is inert in $\Q(\zeta_p)$ for $p = 3,5$. The
second statement follows from the identification of the class
group with the Galois group of the Hilbert class field, and
the fact that all $p$-groups are solvable. $\qed$

\

Assume that $|P| = n = p^k$ for some $k$.
Recall that $M$ is the module generated by $P$ and its conjugates.

\begin{lemma} $\Q(M) = \Q(M[p^k]) \subseteq \Q(\zeta_{p^k})$.
\end{lemma}
By Lemma $1.5$  we find that $M$ is unramified 
outside $p$.
By Theorem $2.2$, $E$ is reducible at $p$.
Since elliptic curves with supersingular reduction at $p$ are
automatically irreducible at $p$, either $E$ has multiplicative
or good ordinary reduction at $p$.
In either case,
the action of inertia at $p$ on $E[p^k]$ is given by
$$(\rho_{p} \mod p^k) |_{\IQ_p} =
\left( \begin{array}{cc} \chi & * \\ 0 & 1 \end{array} \right)$$
where $\chi$ is the cyclotomic character, and $*$ is possibly trivial.
It follows from the almost rationality of $P$ that $M$ is unramified
over $\Q(\zeta_{p^k})$, since:
$$\chi|_{\Gal(\Qbar/\Q(\zeta_{p^k}))} 
 \equiv 1 \mod p^k$$
and so all inertial elements fixing $\Q(\zeta_{p^k})$ will satisfy
$(\sigma - 1)^2 P = 0$. 

\

Since $E$ is reducible at $p$, 
$\Q(E[p],\zeta_{p})/\Q(\zeta_{p})$ is 
a Galois extension of degree dividing $p$.  
$\Gal(\Q(E[p^{k}])/\Q(E[p]))$ is  automatically a $p$-group,
and thus 
$\Gal(\Q(E[p^k])/\Q(\zeta_p))$ is also.  Since 
\mbox{$M
\subseteq E[p^k]$} it follows that the extension
$\Q(M, \zeta_{p^k})/\Q(\zeta_{p^k})$
is an unramified Galois extension whose Galois group is a 
$p$-group.  Thus by Lemma $3.2$, $\Q(M) \subseteq \Q(\zeta_{p^k})$, as claimed.
$\qed$

\

Suppose that $E$ has multiplicative reduction at $p$.  Then locally at
$p$, $E$ is given by a Tate curve $E_q$.  Let $I_p$ denote the
absolute inertia group of $\Q_p$.  For each $n$, we have an exact sequence of
$I_p$-modules:
$$\begin{diagram}
0  \lra \mu_{p^n} & \lra E_q[p^n]  \rTo^{\psi} &
  \Z/p^n\Z \lra  0 \end{diagram}
$$
Let $M' = M \cap \mu_{p^k}$ and $M'' = \psi(M) \subseteq \Z/p^k\Z$.
Then by the Snake Lemma we have an exact sequence of $I_p$-modules: 
$$
0 \lra  M' \lra M \lra M'' \lra 0$$
Since $M$ is defined over $\Q(\zeta_{p^k})$, this is in fact an exact
sequence of $G$-modules, where $G = \Gal(\Q(\zeta_{p^k})/\Q) = \Gal(\Q_p(\zeta_{p^k})/\Q_p) \hookrightarrow \IQ_p$, and where the action of
$\Gal(\Qbar/\Q)$ on $M$ factors through $G$.
Moreover, we observe that as Galois modules, $M''$ is constant and
$M'$ is cyclotomic.

\
   
Now suppose that $E$ has ordinary good reduction at $p$. 
On the level of inertia, multiplicative
and ordinary good reduction are highly analogous.
In particular from
 (\cite{AOEC} VII. Prop. 2.1) and
(\cite{Serre1} Prop. 11), $E[p^n]$ sits in an exact
sequence of $I_p$-modules:
$$\begin{diagram}
0  \lra \mu_{p^n} & \lra E[p^n]  \rTo^{\psi} &
  \Z/p^n\Z \lra  0 \end{diagram}
$$
where $\psi$ is reduction mod $p$ on the elliptic curve. This is
precisely the $I_p$-module sequence arising
when $E$ has multiplicative reduction,  and so we
may treat these cases simultaneously.

\
 
Assume for the moment that $p \geq 3$.  Taking $G$-invariants of the
above sequence we obtain an
exact sequence of Galois cohomology:
$$0 \lra H^0(G,M) \lra H^0(G,M'') \lra H^1(G,M') $$
Since $G$ is cyclic and $p \geq 3$, by Sah's Lemma
(\cite{Sah}, Ch.8 Lemma 8.1), $H^1(G,M') = 0$.
Since $M''$ is a constant module, $H^0(G,M'') = M''$.  Thus we obtain
a map 
$$M'' \simeq H^0(G,M) \hookrightarrow M$$
which induces a
splitting (as Galois modules) of our exact sequence.  Thus $M \simeq
\mu_{p^r} \oplus \Z/p^s\Z$ for some integers $r,s$. From
Lemma $1.2$, it follows that either $P$ is rational, or
$p =3$, $r=1$, and $P$ is a rational point of $3$-power order
plus a point that generates a $\mu_3$ subgroup of $E[3]$. From
Theorem $2.3$, this rational point is of order dividing $9$.

\

Assume now that $p = 2$.  Consider the following sequence of $G$-modules:
$$\begin{diagram}
0  \lra M' & \lra M  \rTo^{\psi} &
  M'' \lra  0 \end{diagram}
$$
Let $H$ be the minimal
quotient of $G$ which acts on $M'$.  Then $H =
\Gal(\Q(\zeta_{2^r})/\Q)$ where $r = \log_2 |M'|$.  We show that the
action of $G$ on $M$ factors through $H$.  Let $\sigma \in G$ become
trivial in $H$ (so $\sigma$ acts trivially on $M'$).
Since $\psi(\sigma P - P) = \sigma \psi(P) - \psi(P) = 0$,
one sees that
$\sigma P - P \in M'$. Thus $(\sigma - 1)^2 P = (\sigma - 1)(\sigma P - P)
  = 0$.
By almost rationality $\sigma P = P$.
Since all the conjugates of $P$ are also almost rational, we see that
$\sigma$ fixes $M$, and the claim follows.  

\

Let $Q \in M$.  The map $H \lra M'$ given by $\sigma \mapsto \sigma Q
- Q$ defines an element of 
$$H^1(H,M') = H^1(\Gal(\Q(\zeta_{2^r})/\Q),
\mu_{2^r}) \simeq H^1((\Z/2^r \Z)^*,\Z/2^r \Z).$$
Sah's Lemma only shows that this group is killed by $2$. In fact, an
elementary argument shows that it is isomorphic to $\Z/2\Z$, where
the non-trivial cocycle class is given by 
$$f(-1 \mod 2^r) = 2^{r-1} \mod 2^r, \qquad f(-3
\mod 2^r) = 0 \mod 2^r.$$
($0$ and $2^{r-1}$ correspond to the elements $1$ and $-1$
in $\mu_{2^r}$).  The cocycle picked out by $Q$ is the image of $Q$
under the composition 
$$\begin{diagram} M & \rTo^\psi & M'' & \rTo^\delta & H^1(H,M') \\
\end{diagram}$$ 
This cocycle is given by $\delta_{\psi(Q)}: \sigma
 \mapsto \sigma Q - Q$.
Since $H^1 \simeq Z/2 \Z$, we may write $\delta_{\psi(Q)}$
 as  $0$ or $f$ 
plus some coboundary.
In particular,
$\sigma Q -  Q = \sigma Q' - Q' +  T_{\sigma}$
for all $\sigma$ and some fixed  (depending only on $Q$) $Q' \in M'$,
and where $T_\sigma$ is
either $0$ or $f(\sigma)$. In particular, $T_{\sigma}$ 
 is Galois invariant (since $\pm 1$ is Galois invariant
in $\mu_{2^r}$) and killed by $2$.  
Thus for any $Q \in M$ we may write 
$Q =Q' + Q''$ with $Q'
\in M'$ such that 
$$\sigma Q = Q + \sigma Q' - Q' + T_{\sigma} = \sigma Q' + Q'' +
T_\sigma.$$
Moreover, since $T_{\sigma}$ 
 is a cocycle, $T_{\tau \sigma} = \tau T_\sigma + T_\tau =
T_\sigma + T_\tau$.

\

Choose $\sigma, \tau \in \Gal(\Qbar/\Q)$ such that:
$$\chi(\sigma)
\equiv 1 + 4k \mod 2^r, \qquad \chi(\tau) \equiv 1 - 4k \mod 2^r$$ 
Note that $(1+4k)(1-4k)
\equiv 1 \mod 8$ and so $\tau \sigma^{-1} = \gamma^2$,
for some $\gamma \in H$.  Thus $T_\tau = T_{\tau \sigma^{-1}} +
T_{\sigma} = 2T_\gamma + T_\sigma = T_\sigma$.  One finds that $\sigma
P + \tau P - 2P = T_\sigma - T_\tau = 0$.  By almost rationality, $P =
\sigma P$. 
Thus the action of $H$ on $P$ factors through
$$\Gal(\Q(\zeta_{2^r})/\Q)/\{\sigma | \chi(\sigma) \equiv 1 \mod 4\}
\simeq \left\{ \begin{array}{ll} \Gal(\Q(i)/\Q), & r \ge 2 \\
1, & r \le 1 \end{array} \right.$$
and so $P$ (and hence $M$) is defined over $\Q(i)$.
In particular, since $H = \Gal(\Q(\zeta_{2^r})/\Q)$ acts faithfully
on $M' \subseteq M$, we conclude that
 $r \le 2$. Since $r \le 2$ and $P' \in M'$, it follows that 
$2 P'$ is either $1$ or $-1$ in $\mu_{2^r}$, and so in particular
is rational.
For any
$\sigma \in H$, 
$$\sigma(2P) = 2 \sigma(P' + P'') = 2 \sigma P' + 2 P'' + 2 T_{\sigma} =
\sigma(2P') + 2P'' = 2P$$
and so $2P$ is rational.  It follows from
Lemma $1.1$(6) that $P$ is rational.

\subsection{Case II: $S = \{2,p\}$, $n = p^m 2^k$, $p = 5,3$} 

From Theorem $2.2$, $E[p]$ is reducible.  We show that $E[2]$ is also
reducible.  Assume otherwise.  From Lemma $1.5$,
$E[2]$ is unramified outside $2$ and $p$, and thus $N(\wtp_2)$ is
equal to $1$ or $p$.  
Since $E[2]$ is irreducible,
Tate's theorem \cite{Tate} eliminates the first possibility.  Thus the
representation $\wtp_2$ is genuinely ramified at $p$, and
from Theorem $2.1$, finite at $2$.  The arguments
of \cite{Ribet1} do not apply at the prime $2$. However, another
result of Ribet \cite{Ribet3} allows us (since $E$ is modular) to
conclude that $\wtp_2$ arises from a modular form of weight $2$ and
level at most $p$.  Since $p \leq 5$ the space of such forms is
trivial, and so $E[2]$ is reducible.

\

Since $E[2]$ is reducible, Corollary $2.1$ shows that $M[p^m]$ is
unramified outside $p$. Arguing as in Lemma $3.3$, we infer
that
$\Q(M[p^m]) \subseteq \Q(\zeta_{p^m})$.
For the proof to go through, it suffices to recall that 
any  $\sigma \in I_p$ satisfying $(\sigma -1)^2 P_p = 0$
 (indeed, any
$\sigma \in I_p$) will automatically satisfy
$(\sigma - 1)^2 P_2 = 0$.

\

\begin{lemma} 
$\Q(M[2^k],\zeta_p)/\Q(\zeta_p)$ is unramified
at all primes above $p$.
\end{lemma}

\Proof Consider the following diagram of
fields:
$$\begin{diagram}
\Q(\zeta_{p^m})	& \rLine &	\Q(\zeta_{p^m},M) & = \Q(\zeta_{p^m},M[2^k]) \\
\dLine		&	&	\dLine & 	\\
\Q(\zeta_p)	& \rLine &	\Q(\zeta_p,M[2^k]) & \\
  \end{diagram}$$
Any element of $\Gal(\Q(M,\zeta_{p^m})/\Q(\zeta_{p^m}))$ (trivially)
fixes
$\Q(\zeta_{p^m})$, and thus $M[p^m]$ 
(since $\Q(M[p^m]) \subseteq \Q(\zeta_{p^m})$). 
Since $E$ is semistable, any $\sigma \in I_p$ satisfies $(\sigma-1)^2 P_2
= 0$. In particular, any inertia at primes above $p$ in
$\Q(M,\zeta_{p^m})/\Q(\zeta_{p^m})$ satisfies
$(\sigma - 1)^2 P_2 = 0$ and $\sigma P_p = 0$, and so
$(\sigma - 1)^2 P = 0$. 
By almost rationality, $\sigma P = P$, and thus the
extension $\Q(\zeta_{p^m},M)/\Q(\zeta_{p^m})$
is unramified at all primes above $p$.
One concludes that that ramification index 
$e_p(\Q(\zeta_{p^m},M)/\Q(\zeta_p))$ is equal to 
$e_p(\Q(\zeta_{p^m})/\Q(\zeta_p)) = p^{m-1}$.
Considering the other part of the diagram, however,
one notes (again using semistability) that the 
ramification of
$\Q(M[2^k])/\Q$ and thus of
$\Q(M[2^k],\zeta_p)/\Q(\zeta_p)$  at $p$ is of $2$-power order.
Since this order must divide $p^{m-1}$, we conclude that
$\Q(\zeta_p,M[2^k])/\Q(\zeta_p)$ is unramified at all primes above $p$. $\qed$

\

From Lemma $1.5$, $M$ is unramified outside $2$ and $p$.
As in the proof of Lemma $3.3$, any inertial elements $\sigma \in I_2$
fixing $\Q(\zeta_{2^k})$ will satisfy $(\sigma-1)^2 P_2 = 0$, and,
by semistability, $(\sigma-1)^2 P_p = 0$. Thus
$\Q(M,\zeta_{2^k})$ 
is unramified at all primes above $2$ in  $\Q(\zeta_{2^k})$.
Combining this with Lemma $3.4$ we infer that
$\Q(M[2^k],\zeta_p,\zeta_{2^k})/\Q(\zeta_p,\zeta_{2^k})$ is
unramified everywhere. Since $\Q(E[2^k])/\Q$ is a $2$-extension,
(as $E[2]$ is reducible)
 we infer from Lemma $3.2$
that
$$\Q(M[2^k]) \subseteq \Q(\zeta_{2^k}, \zeta_p).$$
  Thus we have
shown that $\Q(M) \subseteq \Q(\zeta_{n})$. In particular, the Galois
group $\Gal(\Q(M,\zeta_p)/\Q(\zeta_p))$ decomposes as a product of two
groups, each of which acts trivially on either $M[p^m]$ or
$M[2^k]$. 

\begin{lemma} $P_p$ and $P_2$ are themselves almost rational
over $\Q(\zeta_p)$. \end{lemma}

\Proof  Let
$\sigma, \tau \in \Gal(\Q(\zeta_n)/\Q(\zeta_p))$ satisfy
$\sigma P_p + \tau P_p = 2 P_p$. Then there exist $\sigma',\tau'
\in \Gal(\Q(\zeta_n)/\Q(\zeta_p))$ such that $\sigma = \sigma'$
and $\tau'=\tau$  in
$\Gal(\Q(\zeta_{p^m})/\Q(\zeta_{p}))$, and
$\sigma'=\tau' = 1$ in $\Gal(\Q(\zeta_{2^k},\zeta_p),\Q(\zeta_p))$.
This is because 
$$\Gal(\Q(\zeta_n)/\Q(\zeta_p)) \simeq
\Gal(\Q(\zeta_{p^m})/\Q(\zeta_p)) \oplus \Gal(\Q(\zeta_{2^k},\zeta_p)/
\Q(\zeta_p)).$$
Since  $\Q(\zeta_p,P_2) \subseteq \Q(\zeta_p,\zeta_{2^k})$ and
$\Q(\zeta_p,P_p) \subseteq \Q(\zeta_{p^m})$, it follows that
$\sigma' P_2 = \tau' P_2 = P_2$ and $\sigma' P_p +
\tau' P_p = 2 P_p$ and so $\sigma' P + \tau' P = 2 P$. By almost
rationality, $P = \sigma' P = \tau' P$, 
and so $P_p = \sigma' P_p = \tau' P_p$, and $P_p$ is almost rational
over $\Q(\zeta_p)$. An identical argument works for $P_2$. $\qed$

\

We shall now limit possible $P_p$ as in Case I. Exactly as
in Case I, 
as $G=\Gal(\Q(\zeta_{p^r})/\Q)$ modules, $M[p^m] \simeq
\mu_{p^r}\oplus \Z/p^s\Z$ for some $r,s \leq m$ (this result
only required the fact that $\Q(M[p^m]) \subseteq \Q(\zeta_{p^m})$).

Choose
$\sigma, \tau \in \Gal(\Q(\zeta_{p^m})/\Q(\zeta_p))$ such that
$$\chi(\sigma) \equiv 1 + p \mod p^k, \qquad \chi(\tau) \equiv 1 - p \mod
p^k.$$
We see that $\sigma P_p + \tau P_p - 2 P_p = 0$.  Since $P_p$ is
almost rational over $\Q(\zeta_p)$, this implies that $\sigma P_p =
P_p$, and thus either $P_p$ is $\Q$-rational or $r = 1$ and
 $M[p^m] \simeq \mu_p \oplus
\Z/p^s\Z$.  If $P_p$ is $\Q$-rational then $P_2 = P - P_p$ is almost
rational over $\Q$ and so $P_2$ is also $\Q$-rational, from Case I.  Hence
(possibly after subtracting a rational point of $3$-power order,
allowed by Lemma $1.1$(3)) we
may assume that $P_p$ generates a $\mu_p$ subgroup of $E[p]$.

\

The argument for $n = 2^k$ applies equally well over the field
$\Q(\zeta_p)$ instead of $\Q$, and we may similarly 
conclude that $P_2$ is defined over $\Q(\zeta_p)$.  If $p=3$, then we
are in the $P = Q+R+S$ case of Theorem $1.2$.  Note that such a point is
almost rational, since the only non-trivial Galois conjugate of $P$ is
$\sigma P = -Q+R + \sigma S$, and so $2 \sigma P - 2 P = -Q \neq 0$.
The only thing left to check is that $S$ has order less than
$32$. If $S$ had order divisible by $32$, then $E$ and
the $\Gal(\Qbar/\Q(\sqrt{-3}))$ module generated by $S$
 would give rise to a
a $\Q(\sqrt{-3})$-rational point of $X_0(32)$. 
Since
$X_0(32)$ is explicitly given by the elliptic curve $y^2 = x^3 - x$, we
may  calculate its  $\Q(\sqrt{-3})$-rational points
(equivalently, the $\Q$ rational points of $X_0(32)$ and its
twist $y^2 = x^3 - 9 x$ by $\Q(\sqrt{-3})$). One finds that both
curves have rank zero, and that the torsion points do not correspond to
semistable elliptic curves over $\Q$. The existence of the
$3$-rational point over $\Q(\sqrt{-3})$ suggests that this analysis
could be refined to further limit the order of $S$.

\

It remains to eliminate the possibility that $p = 5$. If $E$ is a semistable
elliptic curve with $\mu_5 \hookrightarrow E[5]$, then there exists a
$5$-isogenous curve $E'$ with a rational point of order $5$. Moreover,
$E[2] \simeq E'[2]$, and so $E'$ also has a rational point of order
$2$. Hence $E'$ has a rational point of order $10$.  Such curves are
parameterized by the genus zero curve $X_1(10)$.  

\

We first  consider the case when $E[2]
\subseteq M[2]$. Since $E[2]$ is reducible, it is defined over some
field of degree $d \le 2$.  Since $\Q(M) \subseteq \Q(\zeta_5)$,
$\Q(E[2])$ must  either be $\Q$ or $\Q(\sqrt{5})$.  Since $E[2] \simeq
E'[2]$, the same must be true of $E'[2]$.  Kubert \cite{Kubert} gives
an explicit parameterization of the genus zero curve $X_1(10)$ in
terms of some uniformizing parameter $f$.  The discriminant 
of the cubic in the Weierstrass equation for $E'$ is a (rational) square
times $(2 f - 1)(4 f^2 - 2 f - 1)$. Hence if $\Q(E'[2]) \subseteq
\Q(\sqrt{5})$ then one of the equations
$$y^2 =(2 f - 1)(4 f^2 -  2 f - 1), \qquad
5 y^2 = (2 f - 1)(4 f^2 -  2 f - 1)$$
must have a rational solution. The first curve is $E_1=20A2$ in
Cremona's tables (\cite{Cremona}), the second its twist $E_2=100A1$.
We find that $E_1(\Q) \simeq \Z/6 \Z$ and $E_2(\Q) \simeq \Z/2 \Z$,
where each torsion point corresponds to a curve such that $\Delta_{E'}
= 0$, and so does not correspond to an actual elliptic curve.  Hence
$E[2] \nsubseteq M[2]$.

\

Suppose now that $E[2]$ does not contain $M[2]$.  If $P$ is of order
$10$, then $5P \in M[2]$ is defined over either $\Q$ or $\Q(\sqrt{5})$.
If $5P$ is not defined over $\Q$, then since $E[2]$ is reducible, it
must contain some non--trivial rational
point and then all of $E[2]$ will be defined over $\Q(\sqrt{5})$.
This situation was considered  above. If $5P \in E(\Q)$, then $\sigma P +
\sigma^2 P = 2 P$, where $\chi_5(\sigma) \equiv 3 \mod 5$ and 
$\chi_5(\sigma^2) \equiv 4 \mod 5$. This contradicts the almost
rationality of $P$. Hence $P$ is of order divisible by $20$. 

Let $Q = mP$ be of order $4$. Then since the Galois module
generated by $Q$ does not contain $E[2]$, it must be cyclic. 
By assumption $E[5]$ also contains a cyclic Galois submodule of
order $5$. Together they generate
a Galois invariant cyclic subgroup of order
$20$, which is impossible, since $X_0(20)$ has no non--cuspidal
rational points. 
 We conclude that $M$ cannot contain a $\mu_5$ subgroup.

\subsection{Case III: Remaining $S$}

For all other possible $n$, we may rule out the existence of $P$ by
the following simple arguments.
Denote by $S$ the set of primes dividing $n$.
\begin{itemize}
\item $S = \{7,5,3,2\},\{7,5,3\},\{7,5,2\},\{7,5\}$.  From Theorem $2.2$,
$E[7]$ is reducible.  Since $(7^2-1,5) = 1$, by Corollary $2.1$,
that $M[5]$ is unramified outside $5$. If $E[5]$ was irreducible,
then it would be finite at $5$ and unramified outside $5$. Arguing
as in the proof of Theorem $2.2$, we infer that $E[5]$ must be
 reducible.   Since $E$ is semistable,
there exists an isogenous curve $E'$ with a rational point of order
$35$.  This contradicts Theorem $2.3$.
\item $S = \{7,3,2\},\{7,3\}$. From Theorem $2.2$, $E[7]$ is reducible.
If $E[3]$ is reducible then there exists an isogenous curve $E'$ with
a rational point of order $21$, which contradicts Theorem $2.3$.
Thus $E[3]$
is irreducible and so finite at $3$ by Theorem $2.2$.  One sees that
$N(\wtp_3)$ is either $1$,$2$,$7$ or $14$.  By level lowering \cite{Ribet1},
the only allowable possibility is that $N(\wtp_3) = 14$.  
In this case, $\wtp_3$
must arise as the Galois representation attached to some modular
form in $S_2(\Gamma_0(14))$ (and trivial character). 
This space is one dimensional, and
corresponds to the elliptic curve $X_0(14)$ of conductor $14$.
Yet from Cremona's tables (\cite{Cremona}) all curves of conductor
$14$
 are reducible at $3$, which implies that $\wtp_3$
must also be reducible.  This is a contradiction.
\item $S = \{5,3,2\},\{5,3\}$. From Theorem $2.2$, $E[5]$ is reducible. If
$E[3]$ is reducible then there exists an isogenous curve $E'$ with a
rational point of order $15$, which contradicts Theorem $2.3$.  Thus $E[3]$ is
irreducible and so finite at $3$ by Theorem $2.1$.  One sees that
$N(\wtp_3)$ is either $1$, $2$, $5$ or $10$.  All cases are
impossible, by Ribet's theorem \cite{Ribet1}.
\item $S = \{7,2\}$.  From Theorem $2.2$,  $E[7]$ is reducible. If
$E[2]$ is reducible then there exists an isogenous curve $E'$ with a
rational point of order $14$, which contradicts Theorem $2.3$.  Thus $E[2]$ is
irreducible and so finite at $2$ by Theorem $2.1$.  One sees that
$N(\wtp_2)$ is either $1$ or $7$.  The first possibility is excluded
by Tate \cite{Tate}.  In the second case, since $\wtp_2$ is genuinely
ramified at $7$, we may use Ribet \cite{Ribet3} to conclude that
$\wtp_2$ arises from some modular form in $S_2(\Gamma_0(14))$.
Again from Cremona's tables (\cite{Cremona}), all elliptic curves
of conductor $14$ are reducible at $2$, and thus $\wtp_2$ is 
also, a contradiction.
\end{itemize}

\end{document}